\documentclass[12pt,a4paper,fleqn]{article}
\usepackage{tikz}
\usepackage{caption}
\usepackage{lscape}

\usepackage{fullpage}
\usepackage{amsmath}
\usepackage{amsthm}
\usepackage{amssymb}
\usepackage{amscd}
\usepackage{plain}
\usepackage{graphicx}
\usepackage{t1enc,epsfig}
\usepackage[mathscr]{eucal}
\usepackage{epstopdf}
\usepackage[german,english]{babel}
\usepackage{stmaryrd}

\usepackage{xcolor}
\usepackage{float}

\usepackage[colorlinks, allcolors=blue]{hyperref}
\usepackage[normalem]{ulem}

\numberwithin{equation}{section}

\begin{document}

\def\diy{\displaystyle}

\def\a{{\alpha}}
\def\b{{\beta}}
\def\Gam{{\Gamma}}
\def\gam{{\gamma}}
\def\del{{\delta}}
\def\eps{{\epsillon}}
\def\veps{{\varepsilon}}
\def\del{{\delta}}
\def\veps {\varepsilon}
\def\vphi{{\varphi}}
\def\vrho{\varrho}
\def\ups{{\upsilon}}
\def\Om {\Omega}

\def\pl {\partial}

\def\bbRd{{\bbRm d}} \def\bbRe{{\bbRm e}} \def\bbRi{{\bbRm i}}
\def\bbRt{{\bbRm t}} \def\bbRv{{\bbRm v}} \def\bbRw{{\bbRm w}}
\def\bbRx{{\bbRm x}} \def\bbRy{{\bbRm y}} \def\bbRz{{\bbRm z}}
\def\bbRD{{\bbRm D}}
\def\bbRO{{\bbRm O}} \def\bbRP{{\bbRm P}} \def\bbRV{{\bbRm V}}

\def\bD{{\mathbf D}}

\def\ovphi{\ov\vphi}
\def\ophi{\ov\phi}

\def\bbA{{\mathbb A}} \def\bbB{{\mathbb B}} \def\bbC{{\mathbb C}}
\def\bbD{{\mathbb D}} \def\bbE{{\mathbb E}} \def\bbF{{\mathbb F}}
\def\bbG{{\mathbb G}} \def\bbH{{\mathbb H}}
\def\bbL{{\mathbb L}} \def\bbM{{\mathbb M}} \def\bbN{{\mathbb N}}
\def\bbP{{\mathbb P}} \def\bbQ{{\mathbb Q}}
\def\bbR{{\mathbb R}} \def\bbS{{\mathbb S}} \def\bbT{{\mathbb T}}
\def\bbU{{\mathbb U}}
\def\bbV{{\mathbb V}} \def\bbW{{\mathbb W}} \def\bbZ{{\mathbb Z}}

\def\bbZd{\bbZ^d}

\def\cA{{\mathcal A}} \def\cB{{\mathcal B}} \def\cC{{\mathcal C}}
\def\cD{{\mathcal D}} \def\cE{{\mathcal E}} \def\cF{{\mathcal F}}
\def\cG{{\mathcal G}} \def\cH{{\mathcal H}} \def\cJ{{\mathcal J}}
\def\cP{{\mathcal P}} \def\cT{{\mathcal T}} \def\cW{{\mathcal W}}
\def\cX{{\mathcal X}} \def\cY{{\mathcal Y}} \def\cZ{{\mathcal Z}}

\def\tA{{\tt A}} \def\tB{{\tt B}}
\def\tC{{\tt C}} \def\tF{{\tt F}} \def\tH{{\tt H}}
\def\tO{{\tt O}}
\def\tP{{\tt P}} \def\tQ{{\tt Q}} \def\tR{{\tt R}}
\def\tS{{\tt S}} \def\tT{{\tt T}}

\def\tx{{\tt x}} \def\ty{{\tt y}}

\def\t0{{\tt 0}} \def\t1{{\tt 1}}

\def\be{{\mathbf e}} \def\bh{{\mathbf h}}
\def\bn{{\mathbf n}} \def\bu{{\mathbf u}}
\def\bx{{\mathbf x}} 
\def\by{{\mathbf y}} 
\def\bz{{\mathbf z}}
\def\bo{{\mathbf o}}
\def\bv{\mathbf v}
\def\ba{{\mathbf a}}
\def\bc{{\mathbf c}}
\def\bb{{\mathbf b}}
\def\B1{{\mathbf 1}} \def\co{\complement}

\def\bmu{{\mbox{\boldmath${\mu}$}}}
\def\bnu{{\mbox{\boldmath${\nu}$}}}
\def\bPhi{{\mbox{\boldmath${\Phi}$}}}

\def\fA{{\mathfrak A}} \def\fB{{\mathfrak B}} \def\fC{{\mathfrak C}}
\def\fD{{\mathfrak D}} \def\fE{{\mathfrak E}} \def\fF{{\mathfrak F}}
\def\fW{{\mathfrak W}} \def\fX{{\mathfrak X}} \def\fY{{\mathfrak Y}}
\def\fZ{{\mathfrak Z}}

\def\bbRA{{\bbRm A}}  \def\bbRB{{\bbRm B}}  \def\bbRC{{\bbRm C}}
\def\bbRF{{\bbRm  F}} \def\bbRM{{\bbRm  M}}
\def\bbRS{{\bbRm S}}
\def\bbRT{{\bbRm T}}  \def\bbRW{{\bbRm W}}

\def\ov{\overline}  \def\wh{\widehat}  \def\wt{\widetilde}

\def\es {{\varnothing}}

\def\bt {\circ}
\def\cc {\circ}

\def\wt{\widetilde}

\def\wtm{{\wt m}} \def\wtn{{\wt n}} \def\wtk{{\wt k}}

\def\be{\begin{equation}}
\def\ee{\end{equation}}

\def\beal{\begin{array}{l}}
\def\beac{\begin{array}{c}}
\def\bear{\begin{array}{r}}
\def\beacl{\begin{array}{cl}}
\def\beall{\begin{array}{ll}}
\def\bealll{\begin{array}{lll}}
\def\beallll{\begin{array}{llll}}
\def\bealllll{\begin{array}{lllll}}
\def\beacr{\begin{array}{cr}}
\def\ena{\end{array}}

\def\bma{\begin{matrix}}
\def\ema{\end{matrix}}

\def\bpma{\begin{pmatrix}}
\def\epma{\end{pmatrix}}

\def\bcs{\begin{cases}}
\def\ecs{\end{cases}}

\def\diy{\displaystyle}

\def\sA{\mathscr A} \def\sB{\mathscr B} \def\sC{\mathscr C}
\def\sD{\mathscr D} \def\sE{\mathscr E}
\def\sF{\mathscr F} \def\sG{\mathscr G} \def\sI{\mathscr I}
\def\sL{\mathscr L} \def\sM{\mathscr M}
\def\sN{\mathscr N} \def\sO{\mathscr O}
\def\sP{\mathscr P} \def\sR{\mathscr R}
\def\sS{\mathscr S} \def\sS{\mathscr T}
\def\sU{\mathscr U} \def\sV{\mathscr V} \def\sW{\mathscr W}
\def\sX{\mathscr X} \def\sY{\mathscr Y} \def\sZ{\mathscr Z}

\def\BZ{{\mathbf Z}}

\def\D{D}

\def\bs {\overline \phi}
\def\bbLv {{\cal E}}

\def\iy{\infty}
\def\ct{\cdot}
\def\cl{\centerline}

\def\bbL{{\mathbb L}} \def\Pf{{\mathbf Z}} \def\f{{\vphi}} \def\g{{\Gam}} \def\s{{\phi}}
\def\boeta{{\mbox{\boldmath$\eta$}}}
\def\sq{\square} \def\tr{\triangle} \def\lz{\lozenge}
\def\cre{\color{red}} \def\cbl{\color{blue}}

\def\bmu{{\mbox{\boldmath${\mu}$}}} 

\def\bbA{{\mathbb A}} \def\bbE{{\mathbb E}} \def\bbH{{\mathbb H}} \def\bbN{{\mathbb N}}
\def\bbR{{\mathbb R}} \def\bbV{{\mathbb V}}  \def\bbW{{\mathbb W}} \def\bbZ{{\mathbb Z}}

\def\cA{{\mathcal A}} \def\cC{{\mathcal C}} \def\cE{{\mathcal E}} \def\cG{{\mathcal G}} 
\def\cN{{\mathcal N}} \def\cP{{\mathcal P}}

\def\Gam{\Gamma} \def\Lam{\Lambda} 
\def\vphi{\varphi} \def\vrho{\varrho}
\def\veps{\varepsilon}
\def\tD{{\tt D}} \def\tR{{\tt R}}
\def\tT{{\tt T}}
\def\wt{\widetilde} \def\ov{\overline}
\def\bbRd{{\bbRm d}} \def\bbRS{{\bbRm S}} 
\def\bbRExt{\bbRm{Ext}}
\def\bbRInt{\bbRm{Int}} \def\bbRSupp{\bbRm{Supp}}


\makeatletter
 \def\fps@figure{htbp}
\makeatother

\def\boldmath{
  \textfont0=\tenbf \scriptfont0=\sevenbf \scriptscriptfont0=\fivebf
  \textfont1=\tenib \scriptfont1=\sevenib \scriptscriptfont1=\fiveib
  \textfont2=\tensyb \scriptfont2=\sevensyb \scriptscriptfont2=\fivesyb
  \textfont3=\tenexb \scriptfont3=\tenexb \scriptscriptfont3=\tenexb
}

\title{{\bf  Minimal Area of a Voronoi Cell in a Packing of Unit Circles }}

\author{ A. Mazel$^1$, I. Stuhl$^2$, Y. Suhov$^{2,3}$}

\date{}

\footnotetext{

 $^1$ \footnotesize{AMC Health, New York, NY, USA;} $^2$ \footnotesize{Math Department, Penn State University, 
PA, USA;} $^3$ \footnotesize{DPMMS, University of Cambridge and St John's College, Cambridge, UK}.
}

\maketitle

\begin{abstract}
We present a new self-contained proof of the well-known fact that the minimal area of a Voronoi cell in a unit circle packing is equal to $2\sqrt{3}$, and the minimum is achieved only on a perfect hexagon. The proof is short and, in our opinion, instructive.
\end{abstract}

\section{Introduction}

This work originated from attempts to rely on the proof in \cite{H} of the fact that the minimal area of a Voronoi cell in a unit circle packing in $\bbR^2$ is equal to $2\sqrt{3}$. Unfortunately, this proof contains some gaps (see Appendix), and our efforts to recover the proof resulted in an alternative approach presented below.

The result itself immediately implies the strong/local version of the theorem of Thue and Fejes T\'oth \cite{Th, To} on dense unit circle packings in $\bbR^2$. Additional applications include the derivation of so-called Peierls estimates in two-dimensional lattice hard-core models of statistical mechanics \cite{MSS1}.

\section{Definitions and basic properties}

For any $\bx \in \bbR^2$ and $r > 0$ denote by $B_r(\bx)$ an open disk of radius $r$ centered at $\bx$. A shorthand notation $B(\bx)$ is used for the {\it unit} disk having $r=1$. A collection of non-overlapping open unit disks is called {\it admissible} and is denoted by $\{B(\bx_i)\}$. An admissible collection $\{B(\bx_i)\}$ represents a {\it unit circle packing}. (It is clear that an admissible collection is finite or countable, as inside each disk there exists a point with rational coordinates.) An admissible collection $\{B(\bx_i)\}$ can be identified with the corresponding collection of centers $\{\bx_i\}$ which is also called admissible. Each element of $\{\bx_i\}$ is called an {\it occupied} point in $\bbR^2$; respectively, we speak about admissible collections of occupied points. Clearly, $|\bx_{i'}-\bx_{i''}| \ge 2$ for any distinct $\bx_{i'}, \bx_{i''} \in \{\bx_i\}$, where $|\bx_{i'}-\bx_{i''}|$ denotes the Euclidean distance between $\bx_{i'},\bx_{i''} \in \bbR^2$.

For each $\bx_{i'} \in \{\bx_i\}$ the corresponding {\it Voronoi cell} is defined as 
$$V(\bx_{i'}) = \left\{\bz \in \bbR^2:\;\;|\bx_{i'}-\bz| \le \inf_{i'' \not = i'}|\bx_{i''}-\bz|\right\}.$$
If for $\bx_{i''}$ the intersection $V(\bx_{i'}) \cap  V(\bx_{i''})$ contains more than one point of $\bbR^2$ then we say that $\bx_{i''}$ is a {\it 
neighboring occupied point} for $\bx_{i'}$ or simply $\bx_{i''}$ is a {\it neighbor} of $\bx_{i'}$. 

Observe that $V(\bx)$ can be unbounded. If $V(\bx)$ is bounded, i.e., $V(\bx) \subset B_{r}(\bx)$ for some $0 < r <\infty$, then $B(\by) \subset B_{4r}(\bx)$ for any neighboring occupied point $\by$. Due to the admissibility requirement the number of such $\by$ cannot exceed $|B_{4r}(\cdot)|/|B(\cdot)| < \infty$, where $|\cdot|$ denotes the area of the corresponding disk.

Our aim is to find the minimal possible area of a Voronoi cell among all Voronoi cells in all admissible collections of occupied points. Suppose that a Voronoi cell with minimal or close to minimal area can be found in an admissible collection containing unbounded Voronoi cells. Then this collection can be completed without breaking the admissibility by a finite or countable set of additional occupied points such that the resulting admissible collection has bounded Voronoi cells only. For that reason from now on we mainly consider admissible collections without unbounded Voronoi cells.

The rest of this section verifies some standard properties of Voronoi cells which makes the entire argument self-contained.

\medskip
{\bf Lemma 1.} {\sl For an occupied point $\bx$ in an admissible collection the corresponding Voronoi cell $V(\bx)$ contains the closure of $B(\bx)$.}

\medskip
{\bf Proof.} If for a point $\bz \in B(\bx)$ there exists an occupied point $\by \not = \bx$ with $|\by-\bz| < |\bx-\bz| < 1$ then by the triangle inequality $|\by-\bx| < 2$ which contradicts the admissibility of the collection. \qed

\medskip
{\bf Lemma 2.} {\sl For an occupied point $\bx$ in an admissible collection the corresponding Voronoi cell $V(\bx)$ is a convex subset of $\bbR^2$.}

\medskip
{\bf Proof.} For an occupied point $\by \not = \bx$ and any $\bz \in V(\bx)$ one has $|\bz-\bx|^2 \le |\bz-\by|^2$ and consequently 
$$|\bx|^2 - 2\bz\cdot \bx \le |\by|^2 - 2\bz\cdot \by.$$ Consider two distinct points $\bz', \bz'' \in V(\bx)$ and suppose that $\bz = \lambda \bz' + (1-\lambda) \bz''$, where $0 \le \lambda \le 1$. Then
$$\beacl
|\bz -\bx|^2 &= |\bz|^2 + \lambda(|\bx|^2-2 \bz'\cdot \bx) + (1-\lambda)(|\bx|^2-2 \bz''\cdot \bx)\cr \\
&\le |\bz|^2 + \lambda(|\by|^2-2 \bz'\cdot \by) + (1-\lambda)(|\by|^2-2 \bz''\cdot \by) \cr
&= |\bz -\by|^2,\ena
$$
which establishes the lemma. \qed

\medskip
{\bf Lemma 3.} {\sl For an occupied point $\bx$ in an admissible collection the boundary of the corresponding Voronoi cell $V(\bx)$ is piecewise linear.}

\medskip
{\bf Proof.} Consider two different occupied points $\bx, \by$ with $V(\bx)\cap V(\by)$ containing more than one $\bbR^2$ point. Let $\bz', \bz'' \in V(\bx)\cap V(\by)$ and $\bz' \not= \bz''$. Now suppose that $\bz = \lambda \bz' + (1-\lambda) \bz''$, $0 \le \lambda \le 1$. Then
$$\beacl
|\bz -\bx|^2 &= |\bz|^2 + \lambda(|\bx|^2-2 \bz'\cdot \bx) + (1-\lambda)(|\bx|^2-2 \bz''\cdot \bx)\cr
&= |\bz|^2 + \lambda(|\by|^2-2 \bz'\cdot \by) + (1-\lambda)(|\by|^2-2 \bz''\cdot \by) \cr
&= |\bz -\by|^2,\ena
$$
i.e., $\bz \in V(\bx)\cap V(\by)$. Also note that for any two distinct occupied points $\bx, \by$ the set
$$\left\{ \bz \in \bbR^2:\;\; |\bz -\bx| \le |\bz -\by| \right\}$$
is a closed half-plane. Consequently, $V(\bx)$ is an intersection of a finite or countable number of closed half-planes. \qed 

\medskip
Lemmas~1-3 are valid for both bounded and unbounded Voronoi cells in an admissible collection of occupied points. If $V(\bx)$ is bounded then these lemmas imply that $V(\bx)$ is a convex polygon containing the unit disk $B(\bx)$. Furthermore, each neighboring occupied point $\by$ is the reflection of $\bx$ with respect to the common side of $V(\bx)$ and $V(\by)$.

The clockwise circular order of sides of a bounded polygon $V(\bx)$ generates the clockwise circular order of the corresponding neighboring occupied points. Connecting the neighboring points in this circular order we obtain a so-called {\it polygon of neighbors}. Together with $\bx$, each side of this polygon uniquely defines a {\it constituting triangle} such that the entire polygon of neighbors is partitioned into constituting triangles. By construction, the vertices of $V(\bx)$ are the centers of the circumcircles of these constituting triangles. The angle of the constituting triangle at vertex $\bx$ is called the {\it constituting angle}.

\medskip
{\bf Lemma 4.} {\sl The circumradius of a constituting triangle is not shorter than $2/\sqrt{3}$.}

\medskip
{\bf Proof.} At least one angle of a constituting triangle, say angle $\alpha$, is not larger than $\pi / 3$. By the admissibility requirement the length $a$ of the opposite triangle side is not shorter than $2$. According to the sine theorem for triangles the triangle circumradius 
$$r = {a \over 2 \sin \alpha} \le {2 \over 2 \sin {\pi \over 3}} = {2\over\sqrt{3}},$$
which establishes the lemma. \qed

\section{Results}

Consider an admissible collection $\{\bx_i\}$ which forms a triangular lattice with the shortest distance between sites equal to  $2$. Then for each $\bx_i$ the corresponding polygon of neighbors is a perfect hexagon with the side length equal to $2$. Correspondingly, $V(\bx_i)$ is a perfect hexagon with the side length equal to $2 / \sqrt{3}$ and the area $|V(\bx_i)|= 2 \sqrt{3}$.

\medskip
{\bf Theorem.} {\sl The minimal area of a Voronoi cell in an admissible configuration of occupied points is equal to $2 \sqrt{3}$. The minimum is achieved only on occupied points $\bx$ having a perfect hexagon with side length $2/\sqrt{3}$ as its Voronoi cell, or equivalently, a perfect hexagon with side length} $2$ as the corresponding polygon of neighbors.

\medskip
{\bf Proof.} In view of properties of a Voronoi cell presented in Lemmas~1-4 the problem is reduced to finding the polygon of minimal area among all polygons $P$ having the following properties:
\begin{description}
\item{(i)} $P$ is convex.
\item{(ii)} $P$ contains a unit disk.
\item{(iii)} The distances from the vertices of $P$ to the center of this disk are not shorter than~${2 \over \sqrt{3}}$.\end{description}

\noindent
The last property is a weaker replacement of the requirement for a Voronoi cell to have the corresponding neighboring occupied points at distances not shorter than  $2$ from each other. It turns out that this weaker requirement is enough to establish the theorem.

For a convex polygon of area $a$ and one of its angles of measure $\alpha$ define the {\it angular density} (with respect to this angle) as the ratio $E={a \over \alpha}$. Now take any polygon satisfying (i)-(iii) and consider several rays originating at the center $\bo$ of the contained unit disk. 
Let the angles between the clockwise consecutive rays be 
smaller than 
$\pi$. These rays partition the polygon into several convex polygons each located inside 
the angle $\alpha_j$ between the corresponding two clockwise consecutive rays. Clearly, the area $a$ of the original polygon $P$ can be calculated as 
$a= \sum_j E_j \alpha_j$, where $E_j$ is the corresponding angular density (with respect to angle $\alpha_j$). 

For any vertex $\bv$ of $P$ 
the convex hull of this vertex and the unit disk centered at $\bo$ belongs to $P$.  
If $|\bo-\bv|=r$ then this convex hull is bounded by two straight segments $[\bv\ba]$ and $[\bv\bc]$ of length $\sqrt{r^2-1}$ and the arc of the unit circle connecting $\ba$ with $\bc$ and having length $2\pi - 2 \arctan\sqrt{r^2-1}$. The area of the quadrilateral $[\bo\ba\bv\bc]$ is equal to $\sqrt{r^2-1}$ and $|\angle \ba\bo\bc|=2\arctan\sqrt{r^2-1}$. Therefore, the corresponding angular density is $E_{\bv} = {\sqrt{r^2-1}\over 2 \arctan\sqrt{r^2-1}}$. Take any point $\bb$ inside the segment $[\bv\bc]$ and consider the angular density $\overline E_{\bv}$ of the smaller quadrilateral $[\bo\ba\bb\bv]$ with respect to the angle $\angle \ba\bo\bb$. Let $|\bb-\bc| = x < \sqrt{r^2-1}$. Then
$$\overline E_{\bv} = {\sqrt{r^2-1} - {x \over 2} \over 2\arctan \sqrt{r^2-1} - \arctan x} > {\sqrt{r^2-1} \over 2\arctan \sqrt{r^2-1}} = E_{\bv}.$$

For two clockwise consecutive polygon vertices $\bv_i$ and $\bv_{i+1}$ consider the corresponding convex hulls and two corresponding quadrilaterals $[\bo\ba_i\bv_i\bc_i]$ and $[\bo\ba_{i+1}\bv_{i+1}\bc_{i+1}]$. Observe that consecutive open quadrilaterals are either adjacent (have a common side) or intersecting.  (In the case when they are separated by 
a non-zero angle $\angle \bc_i\bo\ba_{i+1}$ the segment $[\bv_i \bv_{i+1}]$ intersects the interior of the unit disk which contradicts property (ii).) With this observation at hand, denote by $\bb_i$ the intersection point of segments $[\bv_i\bc_i]$ and $[\bv_{i+1}\ba_{i+1}]$. According to the displayed equation above, the angular density of $[\bo\ba_i\bv_i\bb_i]$ is larger than the angular density of $[\bo\ba_i\bv_i\bc_i]$ and the angular density of $[\bo\bb_i\bv_{i+1}\bc_{i+1}]$ is larger than the angular density of $[\bo\ba_{i+1}\bv_{i+1}\bc_{i+1}]$.

\def\cQ{\mathcal Q}

Consider now the union of the quadrilaterals $[\bo\ba_i\bv_i\bc_i]$ over all vertices $\bv_i$. It is a polygon $Q$ (generally, non-convex) contained in $P$. Between each two vertices $\bv_i$ and $\bv_{i+1}$ the polygon $Q$ may contain an additional vertex $\bb_i$ that is the intersection point introduced above. The angular density of $[\bo\bb_{i-1}\bv_i\bb_i]$ is not smaller than the angular density of $[\bo\ba_i\bv_i\bc_i]$, and the two angular densities are equal only if $\bb_{i-1}=\ba_i$ and $\bb_i=\bc_i$. For $r \ge {2 \over \sqrt{3}}$ the minimum of $E_{\bv} = {\sqrt{r^2-1}\over 2\arctan\sqrt{r^2-1}}$ equals $\sqrt{3} \over \pi$ and is achieved only at $r = {2 \over \sqrt{3}}$. Thus, the total area of the union of the quadrilaterals $[\bo\ba_i\bv_i\bc_i]$ (or equivalently the union of mutually disjoint open quadrilaterals $[\bo\bb_{i-1}\bv_i\bb_{i}]$) is not smaller than ${\sqrt{3} \over \pi} 2\pi = 2\sqrt{3}$. Obviously, this minimum is achieved only when the interiors of the quadrilaterals $[\bo\ba_i\bv_i\bc_i]$ are disjoint (i.e., $[\bo\ba_i\bv_i\bc_i] = [\bo\bb_{i-1}\bv_i\bb_{i}]$) and $|\bo-\bv_i| = {2 \over \sqrt{3}}$ for all $i$. In this case the corresponding polygon is the perfect hexagon with side length ${2 \over \sqrt{3}}$. 

Indeed, if $V(\bx)$ is a hexagon with $r_i > {2 \over \sqrt{3}}$ for some $i$ then $|\angle \ba_i \bo \bc_i| > {2 \pi \over 6}$. Therefore, such a hexagon has the angular density larger than minimal for some non-zero angle. Consequently, its area is larger than $2\sqrt{3}$.

Any polygon with $n > 6$ vertices contains at least two overlapping quadrilaterals $[\bo \ba_i\bv_i\bc_i]$ even if $r_i = {2 \over \sqrt{3}}$ for all $i$ because $n {2 \pi \over 6} > 2\pi$. Therefore, the angular density becomes larger than minimal for some non-zero angle. Hence, any polygon satisfying (i)-(iii) with more than 6 vertices has a non-minimal area. 

A polygon with $n = 3,4$ or 5 vertices necessarily has at least one $|\angle \ba_i \bo \bc_i| \ge {2\pi \over n}$ and therefore $r_i > {1 \over \cos{\pi \over n}} > {2 \over \sqrt{3}}$. Thus, the corresponding angular density $E_{\bv_i}$ is again greater than $\sqrt{3} \over \pi$ for some non-zero angle; consequently, the area of the entire polygon is larger than the minimal area $2\sqrt{3}$.~\qed

\section{Appendix}

The problem with the argument in \cite{H} is that at some point the proof considers the ``most critical case'', but it is not specified what quantity is optimized at this ``most critical case'' and why this quantity is important. The desired quantity seems to be the area excess in $V(\bx)$ over the area of the contained unit disk which can be attributed to a single non-close neighbor. (In the terminology of \cite{H} a neighbor $\by$ is non-close to the center $\bx$ of $V(\bx)$ iff $|\bx-\by| > 2.3$.) By considering this ``most critical case'' for a single non-close neighbor the proof in \cite{H} concludes that the corresponding area excess is at least $0.21$. After that the proof claims that the presence of two non-close neighbors implies that the area excess is at least $0.42$. This additivity assumption is actually wrong as individual excesses can overlap, and one can give a counterexample showing two non-close neighbors with the total attributed area excess smaller than $0.42$.

\end{document}